\documentclass{amsart}
\pdfoutput=1
\usepackage{amsfonts}
\usepackage{amsmath}
\usepackage{amssymb}
\usepackage{graphicx}
\newtheorem{theorem}{Theorem}[section]
\theoremstyle{plain}

\newtheorem{definition}[theorem]{Definition}
\newtheorem{example}[theorem]{Example}

\newtheorem{lemma}[theorem]{Lemma}

\numberwithin{equation}{section}
\begin{document}
\title{On the Necessity of Reidemeister Move 2 for Simplifying Immersed Planar Curves}
\author{Tobias J. Hagge}
\address{Indiana University}
\email{thagge@indiana.edu}
\author{Jonathan T. Yazinski}
\address{Indiana University}
\email{jyazinsk@indiana.edu}

\begin{abstract}
In 2001,\ \"{O}stlund\ conjectured that Reidemeister moves 1 and 3 are
sufficient to describe a homotopy from any generic immersion $S^{1}
\rightarrow\mathbb{R}^{2}$ to the standard embedding of the circle. We show
that this conjecture is false.

\end{abstract}
\maketitle

\section{Introduction}

We wish to consider the problem of simplifying immersed planar curves, in a
sense which will later be made precise. Intuitively, a generic immersion
$S^{1}\rightarrow\mathbb{R}^{2}$ can be considered as a knot diagram without
the crossing data, and for such immersions we can apply planar versions of the
Reidemeister moves for knot diagrams. By applying all three Reidemeister moves
to such a diagram, one is able to obtain a standardly embedded circle with no
double points. One way to see this is to add crossing data so as to give a
knot diagram of the unknot; applying the standard three Reidemeister moves to
this knot diagram gives the standardly embedded circle \cite{[Reid]}.

In \cite{[Oest]}, \"{O}stlund observed that Reidemeister 1 is the only move
that changes the degree of the Gauss map, and showed that Reidemeister move 3
is the only move that can change the signed number of instances of certain
subdiagrams of the Gauss diagram for an embedding. These properties were used
to show that any knot $K$ admits a pair of diagrams such that every sequence
of Reidemeister moves connecting them contains instances of Reidemeister moves
1 and 3. Planar versions of the same arguments give immersions of the circle
in which every connecting sequence contains instances of Reidemeister moves 1
and 3. \"{O}stlund conjectured that any sequence of Reidemeister moves could
be replaced with a sequence consisting of only moves 1 and 3. A counterexample
to this conjecture for the case of knots appears in \cite{[Mant]}.
Independently, the first auther showed in \cite{[Hag]} that every knot type
admits pairs of diagrams such that every connecting sequence contains every
Reidemeister move type.

These two arguments, however, rely heavily on information about the crossings,
and do not generalize to the case of planar Reidemeister moves. The purpose of
the present paper is to disprove the planar version of the conjecture. Since
every sequence of planar Reidemeister moves corresponds to a (non-unique)
sequence of knot Reidemeister moves, the counterexample also provides an
alternate disproof of the conjecture for knots.

\bigskip

\section{Definitions and Main Results}

\begin{definition}
An \emph{immersed curve} is the image of a map $f:S^{1}\rightarrow F,$ where
$F$ is some surface, such that any point in the pair $\left(  F,f\left(
S^{1}\right)  \right)  $ has a neighborhood homeomorphic to a neighborhood in
the picture below. The pair $\left(  F,f\left(  S^{1}\right)  \right)  $ shall
denote the immersed curve.

\begin{center}
\includegraphics[
natheight=1.333500in,
natwidth=1.333500in,
height=0.6953in,
width=0.6953in
]
{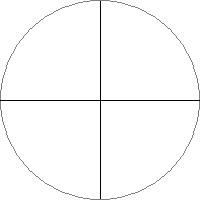}
\\
Figure 1
\end{center}

\end{definition}

In this paper, $F$ will usually be $\mathbb{R}^{2}$ or $S^{2}$.

\begin{definition}
\label{Rmoves}Given an immersed curve, the Reidemeister moves are given, as
numbered below (1a, 1b, 2a, 2b, or 3), by identifying a disk in $\left(
F,f\left(  S^{1}\right)  \right)  $ homeomorphic to the disk on the left side
of the numbered picture and replacing it with the homeomorphic preimage of the
disk on the right.

\begin{center}
\includegraphics[
natheight=7.292100in,
natwidth=16.666700in,
height=1.4857in,
width=3.3615in
]
{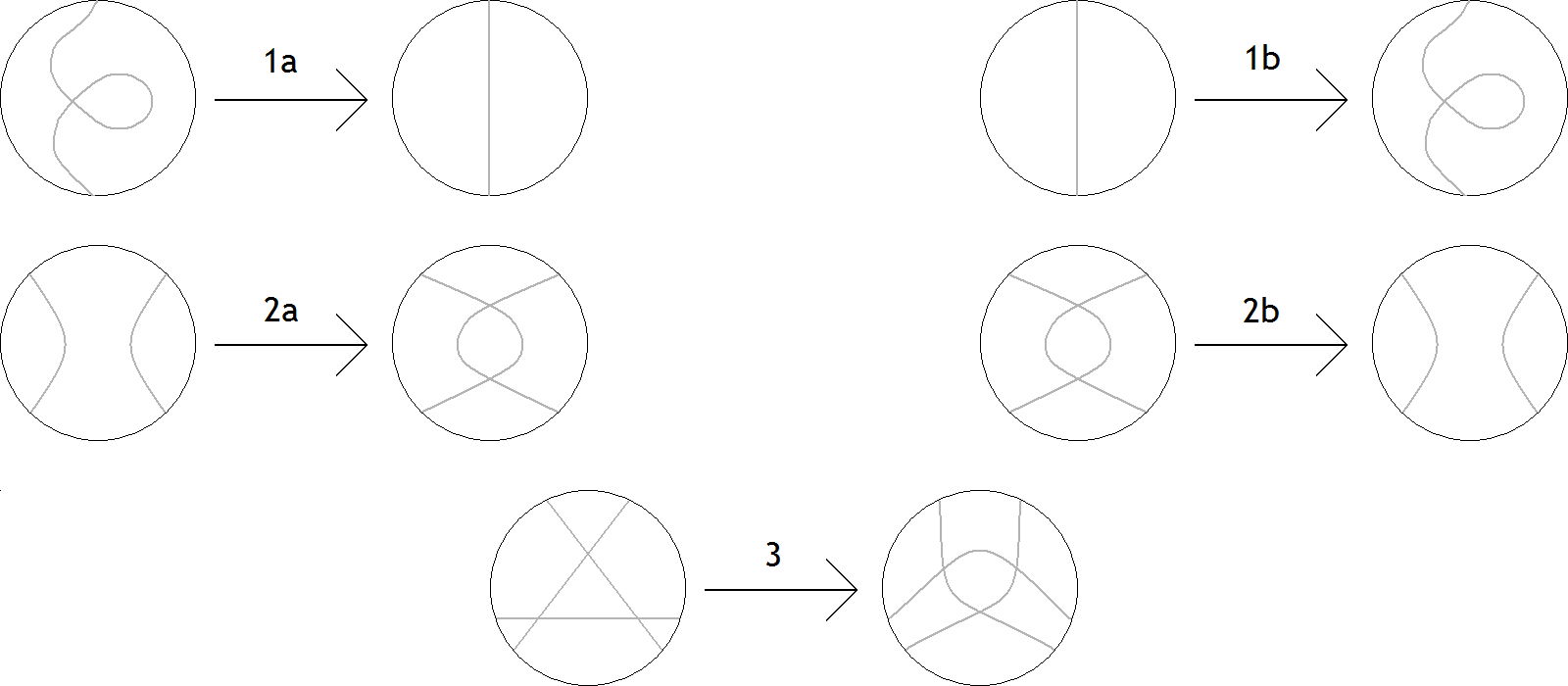}
\\
Figure 2
\end{center}

\end{definition}

By convention, planar isotopies are always allowed as moves, even when not
explicitly mentioned.

\begin{theorem}
Any two homotopic immersed curves are connected by a sequence of Reidemeister
moves and planar isotopies.
\end{theorem}

\begin{proof}
The argument in \cite{[Reid]} works for an arbitrary surface when there is no
crossing data.
\end{proof}

\begin{definition}
An immersed curve $c_{0}$ is \emph{$(1,3)$-simplifiable} if for some $N$ there
exists a sequence of immersed curves $\left\{  c_{i}\right\}  _{i=0}^{N}$ such
that $c_{i+1}$ is obtained from $c_{i}$ by applying one of Reidemeister moves
1a, 1b, or 3, and $c_{N}=\left(  F,f\left(  S^{1}\right)  \right)  ,$ where
$f$ is an embedding. The sequence $\left\{  c_{i}\right\}  _{i=0}^{N}$ is
called a \emph{simplifying sequence} for the curve $c_{0}.$
\end{definition}

\begin{example}
If $F$ is a surface of genus at least $1$ and $c_{0}$ is not nullhomotopic,
then $c_{0}$ is not $(1,3)$-simplifiable. This is because the Reidemeister
moves applied to curves preserve homotopy type.
\end{example}

\"{O}stlund's conjecture, stated in our language, is that every immersed
planar curve is $(1,3)$-simplifiable. Since any curve which is $(1,3)$
-simplifiable in $R^{2}$ is $(1,3)$-simplifiable in its one point
compactification $S^{2}$, the next theorem suffices to disprove the conjecture:

\begin{theorem}
[Main Theorem]\label{Main}The following curve is not $(1,3)$-simplifiable in
$S^{2}$:\newline

\begin{center}
\includegraphics[
natheight=2.866000in,
natwidth=2.866000in,
height=2.034in,
width=2.034in
]
{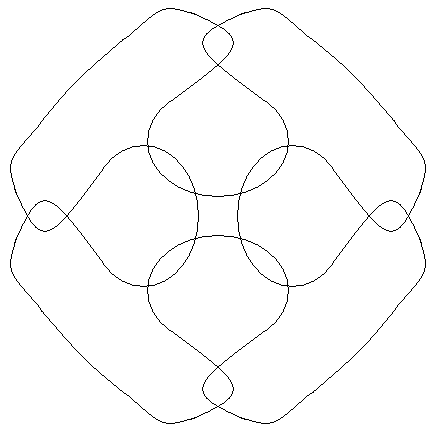}
\\
Figure 3
\end{center}

\end{theorem}

The proof of this theorem does not rely on heavy machinery. It should be noted
that there are immersed curves without an obvious simplifying sequence, which
are nonetheless $(1,3)$-simplifiable. For example, consider the following:

\begin{center}
\includegraphics[
natheight=2.667100in,
natwidth=2.667100in,
height=1.8948in,
width=1.8948in
]
{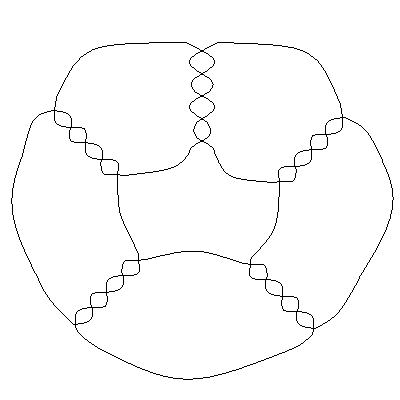}
\\
Figure 4
\end{center}

The easiest way to show that this curve is $(1,3)$-simplifiable is to apply
this theorem:

\begin{theorem}
\label{Double Bigon}Let $c$ be a $(1,3)$-simplifiable curve. Suppose that in
$c$ we replace some instances of the local picture\newline

\begin{center}
\includegraphics[
natheight=1.333500in,
natwidth=1.333500in,
height=0.6953in,
width=0.6953in
]
{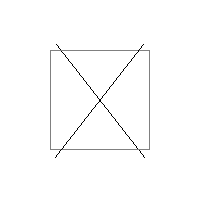}
\end{center}

with the local picture

\begin{center}
\includegraphics[
natheight=1.333500in,
natwidth=1.333500in,
height=0.6953in,
width=0.6953in
]
{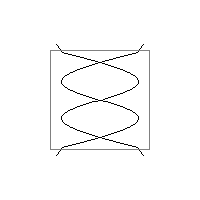}
\end{center}

relative boundary (i.e. double bigons replace double points) to obtain curve
$c^{\prime}$. If $c$ is $(1,3)$-simplifiable, then $c^{\prime}$ is $(1,3)$-simplifiable.
\end{theorem}

It should be noted that Theorem~\ref{Double Bigon} does not say that moves 1
and 3 may be used to replace a double point with a double bigon in an
arbitrary diagram. Nonetheless, applying Theorem~\ref{Double Bigon} repeatedly
to Figure 5 gives the following $(1,3)$-simplifiable immersed curve:

\begin{center}
\includegraphics[
natheight=2.000300in,
natwidth=2.000300in,
height=1.0274in,
width=1.0274in
]
{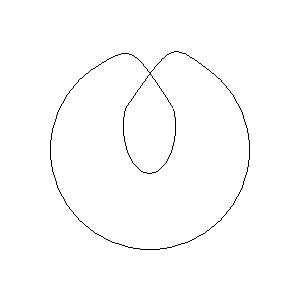}
\end{center}

One could generalize \"{O}stlund's conjecture and ask whether two homotopic
curves on a surface are related by only the first and third Reidemeister
moves. This generalized conjecture is much easier to falsify. It is in fact a
generalization because all generic curves on $\mathbb{R}^{2}$ or $S^{2}$ are
homotopically trivial.

\begin{theorem}
\label{Torus}The following two curves on $T^{2}$ are homotopic, but are not
related by a sequence of Reidemeister moves consisting of only the first and
third moves. \newline

\begin{center}
\includegraphics[
natheight=1.333500in,
natwidth=3.333000in,
height=1.3681in,
width=3.3797in
]
{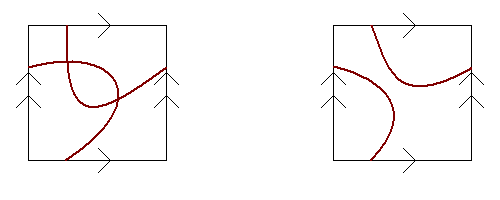}
\end{center}

\end{theorem}

\section{Proof of Main Theorem}

This section proves Theorem \ref{Main}. Consider the following shaded regions
in the curve from Figure~4, interpreted as a diagram on $S^{2}$:

\begin{center}
\includegraphics[
natheight=4.478900in,
natwidth=4.478900in,
height=2.2675in,
width=2.2675in
]
{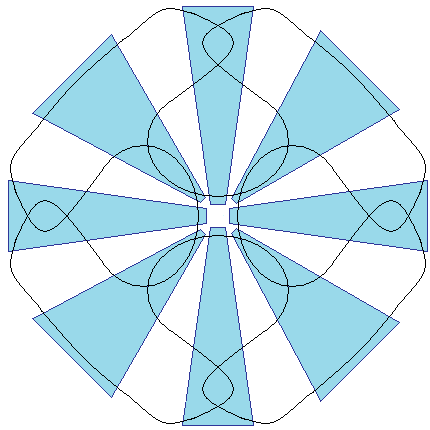}
\\
Figure 5
\end{center}

Reinterpret the diagram as a collection of eight blue boxes containing
immersed tangles, connected by lines with no double points. Each box has a
left and right side, as labeled below; the left side of a given box is
connected to the right side of its neighbor. Two polygons in the diagram
deserve special attention and are marked with a star.

\begin{center}
\includegraphics[
natheight=4.478900in,
natwidth=4.478900in,
height=2.2675in,
width=2.2675in
]
{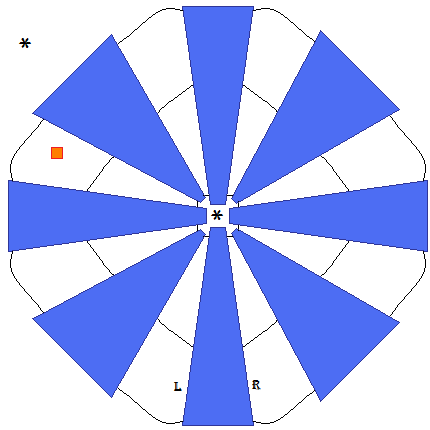}
\\
Figure 6
\end{center}

The diagram satisfies the following properties:

\begin{enumerate}
\item \label{strands ends} Each blue box contains a tangle with three strands.
One of the strands, denoted strand 1, begins and ends at the left side. Strand
2 begins and ends at the right side. Strand 3 has one endpoint on each side of
the box.

\item \label{strands connect} In each box, the left side of strand 3 connects
to strand 2 in the adjacent box to the left. The right side of strand 3
connects to strand 1 in the adjacent box on the right.

\item \label{strands two double} Strands 1 and 2 intersect in exactly two
double points.

\item \label{strands four sides}The polygons marked with a star have at least
four sides.
\end{enumerate}

We will show that any application of moves 1a, 1b, or 3 to any copy of
Figure~6 with immersed tangles satisfying the above properties results in a
diagram which may be interpreted as a copy of Figure~6 with immersed tangles
still satisfying those properties. Since property~\ref{strands two double}
implies that any blue box has at least two double points, every sequence of
such moves results in a diagram with at least sixteen double points. This
proves that the curve is not $(1,3)$-simplifiable.

First, note that a move of type 1a, 1b, or 3 occuring entirely within one of
the blue boxes gives a diagram (with the same boxes) satisfying all of the
above properties. Such a move fixes the endpoints of the strands, so
Properties~\ref{strands ends}~and~\ref{strands connect} remain satisfied. None
of these moves change the number of times one strand intersects another within
a box, so Property~\ref{strands two double} holds after a move.

Property~\ref{strands four sides} actually follows from the arrangement of the
boxes and the other three properties. Fix a starred polygon and consider the
portion of its boundary lying within a single blue box. If the end points of
that boundary portion belong to different strands within the box, that box
contains at least one vertex for the starred polygon. Otherwise,
Property~\ref{strands ends} implies that both ends belong to strand 3. Then
Property~\ref{strands connect} implies that the end points for the portion of
the starred polygon's boundary lying in each of the adjacent blue boxes belong
to different strands within that box. Thus each of the adjacent boxes contains
a vertex for the starred polygon. Therefore, allowed moves cannot reduce the
number of vertices (or edges) for a starred polygon below four.

It remains to show that it suffices to consider only moves lying within a
single blue box. First, consider Reidemeister move 1b. Performing this move
requires a disk in our immersed curve that is homeomorphic to the disk on the
left side of picture 1b in Figure 2 in Definition \ref{Rmoves}. Suppose that
the segment on the left side of picture 1b in Definition \ref{Rmoves} is not
contained completely inside one of the blue boxes as specified above. Then one
can redefine the blue box before performing the move so that it occurs
entirely within a single blue box. For example, suppose that the disk for move
1b is the following:

\begin{center}
\includegraphics[
natheight=4.166700in,
natwidth=4.166700in,
height=1.2782in,
width=1.2782in
]
{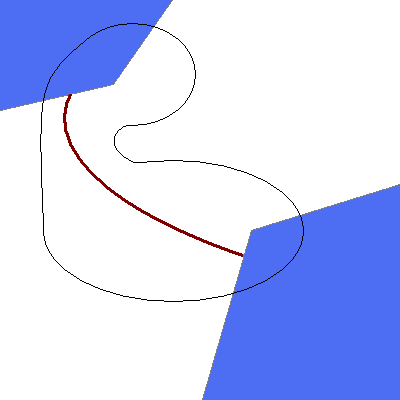}
\end{center}

One can then isotop the blue boxes, while leaving crossings fixed, as follows:

\begin{center}
\includegraphics[
natheight=4.166700in,
natwidth=4.166700in,
height=1.2782in,
width=1.2782in
]
{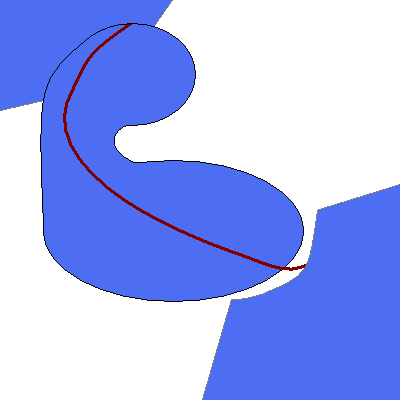}
\end{center}

Move 1a, on the other hand, removes a one-sided polygon. This polygon must lie
entirely within a single blue box, for the following reason: Clearly, a
one-sided polygon cannot separate the two starred regions. If a closed smooth
subcurve of $f\left(  S^{1}\right)  $ does not lie in a single blue box, and
does not separate the two starred regions, then it must enter and exit one of
the blue boxes on the same side. Such a curve contains a segment of strand
type 1 or 2, and by Property \ref{strands two double}, any such curve will
have at least two crossings.

Finally, move 3 always occurs on a neighborhood of a triangle (which can never
be marked with a star). If that triangle lies entirely in one blue box, that
box may be isotoped as above to include the entire disk on which the move
occurs. Otherwise, the triangle intersects the white region (an example of
such a potential triangle is marked with a red dot in Figure~6). One can
verify that this implies that one of the blue boxes intersects the triangle
only in a single corner, as shown in Figure~7, for example.

\begin{center}
\includegraphics[
natheight=3.125400in,
natwidth=3.125400in,
height=1.2782in,
width=1.2782in
]
{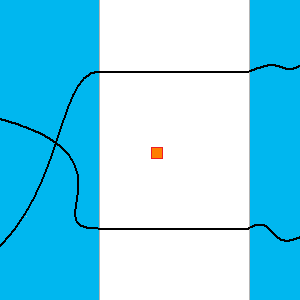}
\\
Figure 7
\end{center}

Assume without loss of generality that the triangle extends to the right of
the blue box containing just the corner, as in Figure 7. There are two
possibilities for the strand ends on the right side of the leftmost box shown
in Figure~7. Either exactly one of the ends belongs to strand 3, or both ends
belong to strand 2. In the box to the right, either both of the pictured left
ends belong to strand 1, or exactly one belongs to strand 3, respectively. In
either case, isotoping the blue boxes in Figure~7 to give the blue boxes in
Figure~8 preserves the required properties and reduces the number of white
regions in the triangle. After at most two such box adjustments, all three
vertices of the triangle must lie in the same blue box. Then, since there are
no isolated blue corners, the entire triangle must lie within a single blue box.

\begin{center}
\includegraphics[
natheight=3.125400in,
natwidth=3.125400in,
height=1.2782in,
width=1.2782in
]
{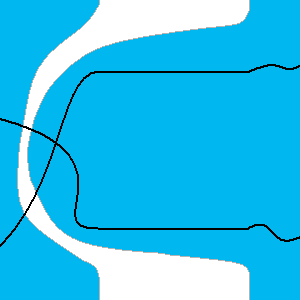}
\\
Figure 8
\end{center}

One could also prove this theorem using Gauss diagrams. We give the main
outline, leaving the proof to the reader. The Gauss diagram for the immersed
curve in Figure 4 is as follows:

\begin{center}
\includegraphics[
natheight=5.562500in,
natwidth=5.791700in,
height=2.2528in,
width=2.3445in
]
{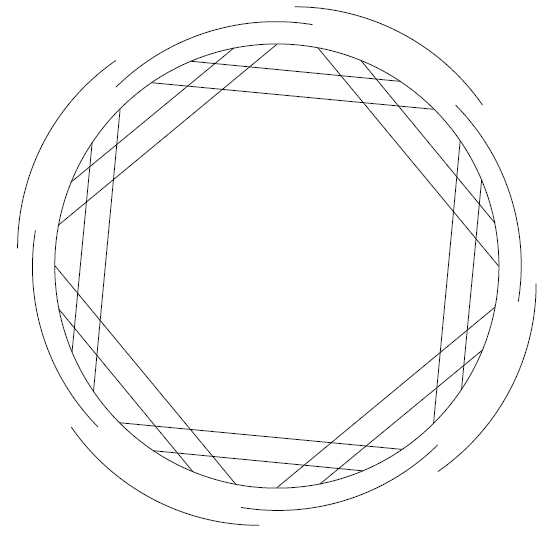}
\\
Figure 9
\end{center}

Move 1b adds a chord to the diagram, which by convention shall be colored
gray. The following properties are preserved by Reidemeister moves 1 and 3.

\begin{enumerate}
\item If all the gray lines are erased, the resulting diagram is exactly as
shown in Figure 9 above except that some of the adjacent and parallel pairs of
black lines may be replaced with crossed pairs.

\item Both endpoints of every gray chord lie in one of the eight regions
indicated in Figure~9.
\end{enumerate}

\section{Proof of Other Theorems\label{Proof Section}}

\begin{lemma}
\label{Double Bigon Lemma}The following pictures are connected by a sequence
of Reidemeister moves 1 and 3: \newline
\begin{center}
\includegraphics[
natheight=1.333500in,
natwidth=3.333000in,
height=0.6953in,
width=1.6942in
]
{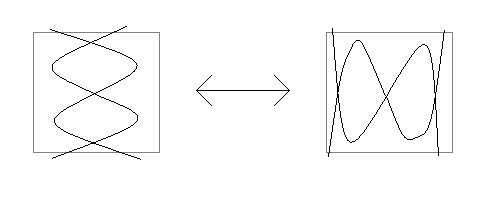}
\end{center}

\end{lemma}

\begin{proof}
This is the necessary sequence of Reidemeister moves:

\begin{center}
\includegraphics[
natheight=3.999800in,
natwidth=2.667100in,
height=2.8262in,
width=1.8948in
]
{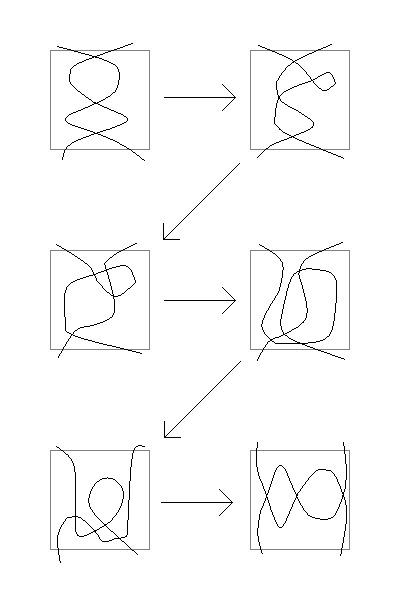}
\end{center}

\end{proof}

\begin{proof}
[Proof of Theorem \ref{Double Bigon}]Consider two immersed curves $L$ and $R$,
equal except inside of a box. The contents of the box for the curves $L$ and
$R$ are given respectively by the following pictures on the left and right:

\begin{center}
\includegraphics[
natheight=1.333500in,
natwidth=3.999800in,
height=0.6953in,
width=2.0271in
]
{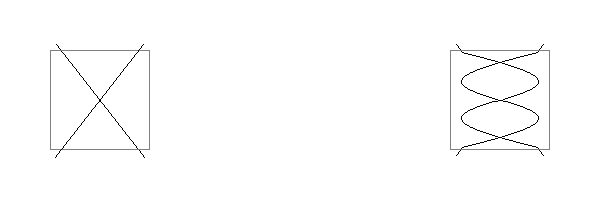}
\end{center}

Suppose $L$ is $(1,3)$-simplifiable. Then Reidemeister moves performed on $L$
that are supported away from the box have analagous Reidemeister moves on $R$.
However, the simplifying sequence for $L$ may contain moves 1a and 3 which
involve the box. The following sequences of moves on $R$ are analogous to
moves 1a and 3 on $L$ which involve the box. In these sequences it may be
necessary to first apply Lemma \ref{Double Bigon Lemma} to obtain the leftmost picture.

\begin{center}
\includegraphics[
natheight=2.083300in,
natwidth=8.333300in,
height=1.0706in,
width=4.1943in
]
{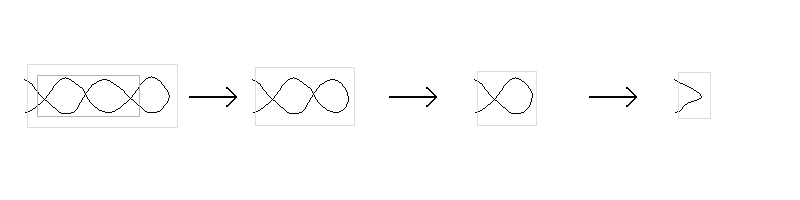}
\end{center}

\begin{center}
\includegraphics[
natheight=2.083300in,
natwidth=8.333300in,
height=1.0706in,
width=4.1943in
]
{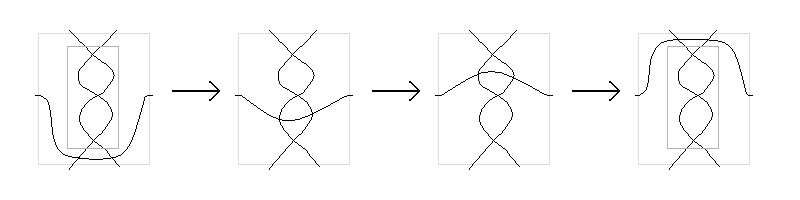}
\end{center}

Applying the moves on $R$ analogous to the moves in a simplifying sequence for
$L$ gives a simplifying sequence for $R$.
\end{proof}

\begin{proof}
[Proof of Theorem \ref{Torus}]It will be sufficient to show that by applying
R1 and R3 moves to a curve on $T^{2}$ of the form

\begin{center}
\includegraphics[
natheight=1.333500in,
natwidth=1.333500in,
height=1.3681in,
width=1.3681in
]
{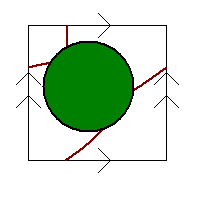}
\end{center}

such that the two strands of the tangle inside the disk intersect, one can
never obtain an embedded curve (i.e. a curve without double points). Observe
that in the picture above, there is exactly one region not contained in the
disk, and this region has genus. Up to isotopy every R1a, R1b and R3 move is
supported within the disk. However, as noted in the proof of the main theorem,
R1 and R3 moves on tangles do not change the number of intersections between strands.
\end{proof}

\end{document}